\documentclass[a4paper,12pt]{amsart}
\usepackage{amsmath}
\usepackage{amsthm}
\usepackage{amscd}
\usepackage{amssymb}
\usepackage{amsfonts}
\usepackage{latexsym}
\usepackage[mathscr]{eucal}

\newtheorem{thm}{Theorem}[section]
\newtheorem{prop}[thm]{Proposition}
\newtheorem{lemma}[thm]{Lemma}

\theoremstyle{definition}

\theoremstyle{remark}

\numberwithin{equation}{section}

\def\supp{\mathop{\rm {supp}}\nolimits}

\def\length{\mathop{\rm {length}}\nolimits}

\def\supp{\hbox{\rm supp}}
\def\supp{\mathop{\rm {supp}}\nolimits}

\def\length{\mathop{\rm {length}}\nolimits}

\def\supp{\hbox{\rm supp}}

\def\Ker{\hbox{\rm Ker}}

\newcommand{\cC}{{\mathcal C}}

\newcommand{\cG}{{\mathcal G}}
\newcommand{\cO}{{\mathcal O}}
\newcommand{\cP}{{\mathcal P}}

\begin{document}
\title{$C^*$-algebras associated with algebraic correspondences 
on the Riemann sphere}

\author{Tsuyoshi Kajiwara}
\address[Tsuyoshi Kajiwara]{Department of Environmental and 
Mathematical Sciences, 
Okayama University, Tsushima, 700-8530,  Japan}      

\author{Yasuo Watatani}
\address[Yasuo Watatani]{Department of Mathematical Sciences, 
Kyushu University,
Hakozaki, Fukuoka, 812-8581,  Japan}

\maketitle
\begin{abstract} Let $p(z,w)$ be a  polynomial in two variables. We 
call the solution of the algebraic equation $p(z,w) = 0$ the 
algebraic correspondence. We regard it as the graph of the 
multivalued function $z \mapsto w$ defined implicitly  by 
$p(z,w) = 0$. 
Algebraic correspondences on the Riemann sphere $\hat{\mathbb C}$ 
give a generalization of dynamical systems of Klein groups and 
rational functions. We introduce $C^*$-algebras associated 
with algebraic correspondences on the Riemann sphere. We show 
that if an algebraic correspondence is free and expansive on a closed 
$p$-invariant subset $J$ of $\hat{\mathbb C}$, then the 
associated $C^*$-algebra 
${\mathcal O}_p(J)$ is simple and purely infinite.

\end{abstract}

\section{Introduction}
For a branched covering $\pi : M \rightarrow M$, 
Deaconu and Muhly \cite{DM} introduced a $C^*$-algebra 
$C^*(M,\pi)$ as the $C^*$-algebra of the r-discrete 
groupoid constructed by Renault \cite{R}. 
In order to capture information of the branched points for the complex 
dynamical system arising from a rational function $R$,  
in \cite{KW1} we introduced a slightly different 
$C^*$-algebras 
${\mathcal O}_R(\hat{\mathbb C})$,  
${\mathcal O}_R = {\mathcal O}_R(J_R)$ and ${\mathcal O}_R(F_R)$
associated with a rational function $R$ on the Riemann sphere,  
the Julia set $J_R$ and the Fatou set $F_R$ of $R$. We showed that 
the $C^*$-algebfras ${\mathcal O}_R(J_R)$ on the 
Julia set is always simple and purely infinite if 
the degree of $R$ is at least two. We also studied a relation 
between branched points and KMS states in \cite{IKW}. 
 C. Delaroche \cite{A} and M. Laca - J. Spielberg \cite{LS} 
showed that a certain boundary action of a Kleinian group 
on the limit set yields a simple nuclear purely infinite 
$C^*$-algebra as groupoid $C^*$-algebra or crossed product. 
Dutkay and Jorgensen study a spectral theory on Hilbert spaces 
built on general finite-to-one maps (\cite{DJ}).
 
On the other hand, Sullivan discovered a dictionary between the theory of 
complex analytic iteration and the theory of Kleinian groups in \cite{Su}. 
Sullivan's dictionary shows a strong analogy between the limit set 
$\Lambda_{\Gamma}$ of a Kleinian group $\Gamma$ 
and the Julia set $J_R$ of a rational function $R$.  Therefore it is natural 
to generalize both Kleinian groups and rational maps.  In fact there 
exist such objects called algebraic correspondences or holomorphic 
correspondences in many works by  several people, for example, 
Bullet \cite{Bu}, Bullet-Penrose \cite{BP1}, \cite{BP2} and 
M\"{u}nzner-Rasch  \cite{MR}. Let $p(z,w)$ be a  polynomial 
in two variables. Then the solution of the algebraic equation 
$p(z,w) = 0$ is called the algebraic correspondence. 
We regard it as the graph of the 
multivalued function $z \mapsto w$ defined implicitly  by 
$p(z,w) = 0$. 

In this paper, 
we introduce $C^*$-algebras associated 
with algebraic correspondences on the Riemann sphere. We show 
that if an algebraic correspondence is free and expansive on a closed 
$p$-invariant subset $J$ of $\hat{\mathbb C}$, then the 
associated $C^*$-algebra 
${\mathcal O}_p(J)$ is simple and purely infinite. We shall 
show some examples and compute the K-groups of the associated 
$C^*$-algebras. For example, 
Let 
$p(z,w) = (w - z^{m_1})(w -  z^{m_2})\dots (w - z^{m_r}) $ 
and $m_1, \dots, m_r$ are all different, where $r$ is the 
number of irreducible components. 
Then $J := {\mathbb T}$ is a $p$-invariant set. 
Let $b = \ ^{\#}B(p)$ be the number of the branched poins. 
Then we have 
\[
 K_0({\mathcal O}_p({\mathbb T})) = {\mathbb Z}^{b}, \ \ 
\text{ and } \ \   
 K_1({\mathcal O}_p({\mathbb T})) = {\mathbb Z}/(r-1){\mathbb Z}.
\]
If $m_1, m_2, \dots, m_r$ are relatively prime, 
then the associated $C^*$-algebra ${\mathcal O}_p({\mathbb T})$ 
is simple and purely infinite. 

Our $C^*$-algebras 
${\mathcal O}_p(J)$ are related with 
$C^*$-algebras of irreversible dynamical systems by 
Exel-Vershik \cite{EV}, $C^*$-algebras associated with subshifts 
by Matsumoto \cite{Ma}, 
graph $C^*$-algebras \cite{KPRR} and their generalizaion for  
topological relations by Brenken \cite{Br}, topological graphs by 
Katsura \cite{Ka1}, and topological quivers by 
Muhly and Solel \cite{MS} and by Muhly and Tomforde \cite{MT}. 
Some of our $C^*$-algebras are isomorphic to $C^*$-algebras 
associated self-similar sets \cite{KW2} and 
Mauldin-Williams graphs \cite{IW}.

\section{Algebraic correspondences}

Let $p(z,w) \in {\mathbb C}[z,w]$ be a  polynomial in two variables  
of degree $m$ in $z$ and degree $n$ in $w$. We shall study  an  
algebraic function implicitly determined by 
the algebraic equation $p(z,w)=0$ on the 
Riemann sphere $\hat{\mathbb C}$. Note that 
there exist two different ways to compactify the algebriac curve 
$p(z,w)=0$. The standard construction in algebraic geometry is 
to consider the zeroes of a homogeneous polynomial $P(z,w,u)$ 
in the complex projective plane ${\mathbb C}P^2$. But we 
choose the second way after \cite{BP1} and introduce four 
variables $z_1, z_2, w_1, w_2$ and a polynomial 
\[
\tilde{p}(z_1,z_2, w_1,w_2) =
z_2^mw_2^np(\frac{z_1}{z_2},\frac{w_1}{w_2}), 
\]
which is separately homogeneous in $z_1,z_2$  and 
in $w_1,w_2$. We identiy  the 
Riemann sphere $\hat{\mathbb C}$ with 
the  complex projective line ${\mathbb C}P^1$.
We denote by $[z_1,z_2]$ an element of ${\mathbb C}P^1$. 
Then the {\it algebraic correspondence} $\cC_p$ of 
$p(z,w)$ on the 
Riemann sphere is a 
closed  subset of $\hat{\mathbb C} \times \hat{\mathbb C}$ 
defined by    
\[
\cC_p := 
\{([z_1,z_2], [w_1,w_2]) \in \hat{\mathbb C} 
\times \hat{\mathbb C} \ | \ 
\tilde{p}(z_1,z_2,w_1,w_2)=0\}. 
\]
Then  $\cC_p$ is compact. In fact,  it is a continuous image of 
a compact subset 
\[
\{(z_1,z_2,w_1,w_2) \in {\mathbb C}^4 \ | \  
\tilde{p}(z_1,z_2,w_1,w_2)=0 \text{ and } 
|z_1|^2 + |z_2|^2 + |w_1|^2 + |w_2|^2 = 1 \}.
\]

To simplify notation, we  write it by 
\[
\cC_p = \{(z,w) \in \hat{\mathbb C} \times \hat{\mathbb C} \ | \ p(z,w) = 0 \}
\]
for short if no confusion can arise. It is also convienient to consider 
change of variables  $u = \frac{1}{z}$  or  $v = \frac{1}{w}$ instead.

For example, let $R(z) = \frac{P(z)}{Q(z)}$ be a rational function with 
polynomials $P(z)$, $Q(z)$. Put $p(z,w)=Q(z)w-P(z)$. 
Then the algebraic correspondence $\cC_p$ of 
$p(z,w)$ on the Riemann sphere is exactly the graph 
$\{(z,w) \in \hat{\mathbb C} \times \hat{\mathbb C}
 \ | \ w = R(z), \ z \in \hat{\mathbb C}\}$ of $R$. 

Therefore 
we regard  the algebraic correspondence $\cC_p$ of a general polynomial 
$p(z,w)$ as the graph of the algebraic function $z \mapsto w$ 
implicitly defined  by the equation 
$p(z,w) = 0$. Then the iteration of the algebraic function is described 
naturally by a sequence $z_1,z_2, z_3, \dots $ satisfying 
$(z_k,z_{k+1}) \in \cC_p$ for $k = 1,2,3, \dots$. 

Any non-zero polynomial $p(z,w) \in {\mathbb C}[z,w]$ has a unique 
factorization into irreducible polynomials:
\[
 p(z,w) = g_1(z,w)^{n_1}\cdots g_p(z,w)^{n_p}
\]
where each $g_i(z,w)$ is irreducible and $g_i$ and $g_j \ (i \not= j)$ 
are prime each other.  

Throughout the paper, we assume that any polynomial $p(z,w)$ we consider is reduced, that is, the above powers $n_i=1$  for any  $i$. We also assume that 
any  $g_i(z,w)$  is not a function only in $z$  or $w$. In particular 
the degree $m$ in $z$ and the degree $n$ in $w$ of $p(z,w)$ are both 
greater than or equal to one.

We need to recall an elementary fact as follows:

\medskip\par\noindent
{\bf Definition}@(branch index). 
Let $p(z,w)$ be a  non-zero polynomial in two variables 
of degree $m$ in $z$ and degree $n$ in $w$. Then we sometimes 
rewrite $p(z,w)$ as 

\begin{align*}
p(z,w) 
 & =   a_m(w)z^m + a_{m-1}(w)z^{m-1}+ \cdots + a_1(w)z + a_0(w) \\
 & = b_n(z)w^n + b_{n-1}(z)w^{n-1}+ \cdots + b_1(z)w + b_0(z) .\\
\end{align*}

 Fix   
$w=w_0 \in \hat{\mathbb C}$. Then the equation $f(z) :=p(z,w_0)=0$ 
in $z \in \hat{\mathbb C}$ has $m$ roots with multiplicities. 
Take any root $z =z_0$. 
The branch index of $p(z,w)$ at $(z_0,w_0)$, denoted by 
$e_p(z_0,w_0)$ or $e(z_0,w_0)$, is defined to be the 
multiplicity for the root $z =z_0$ of $f(z) =p(z,w_0)=0$.  
For example, let $R(z) = \frac{P(z)}{Q(z)}$ be a rational function with 
polynomials $P(z)$, $Q(z)$. Put $p(z,w)=Q(z)w-P(z)$. Then the branch index 
$e_p(z_0,R(z_0))$ coincides with 
the usual branch index $e_R(z_0)$ of $R$ at $z = z_0$.

\section{associated $C^*$-algebras}

We recall Cuntz-Pimsner algebras \cite{Pi}.  
Let $A$ be a $C^*$-algebra 
and $X$ be a Hilbert right $A$-module.  We denote by $L(X)$ be 
the algebra of the adjointable bounded operators on $X$.  For 
$\xi$, $\eta \in X$, the "rank one" operator $\theta _{\xi,\eta}$
is defined by $\theta _{\xi,\eta}(\zeta) = \xi(\eta|\zeta)$
for $\zeta \in X$. The closure of the linear span of rank one 
operators is denoted by $K(X)$.   We say that 
$X$ is a Hilbert bimodule (or $C^*$-correspondence) over $A$ 
if $X$ is a Hilbert right  $A$-
module with a homomorphism $\phi : A \rightarrow L(X)$.  
In this note, we assume 
that $X$ is full and $\phi$ is injective.

   Let $F(X) = \oplus _{n=0}^{\infty} X^{\otimes n}$
be the Fock module of $X$ with a convention 
$X^{\otimes 0} = A$. 
 For $\xi \in X$, the creation operator
$T_{\xi} \in L(F(X))$ is defined by 
$$
T_{\xi}(a) =  \xi a  \qquad \text{and } \ 
T_{\xi}(\xi _1 \otimes \dots \otimes \xi _n) = \xi \otimes 
\xi _1 \otimes \dots \otimes \xi _n .
$$
We define $i_{F(X)}: A \rightarrow L(F(X))$ by 
$$
i_{F(X)}(a)(b) = ab \qquad \text{and } \ 
i_{F(X)}(a)(\xi _1 \otimes \dots \otimes \xi _n) = \phi (a)
\xi _1 \otimes \dots \otimes \xi _n 
$$
for $a,b \in A$.  The Cuntz-Toeplitz algebra ${\mathcal T}_X$ 
is the $C^*$-subalgebra of $L(F(X))$ generated by $i_{F(X)}(a)$
with $a \in A$ and $T_{\xi}$ with $\xi \in X$.  
Let $j_K : K(X) \rightarrow {\mathcal T}_X$ be the homomorphism 
defined by $j_K(\theta _{\xi,\eta}) = T_{\xi}T_{\eta}^*$. 
We consider the ideal $I_X := \phi ^{-1}(K(X))$ of $A$. 
Let ${\mathcal J}_X$ be the ideal of ${\mathcal T}_X$ generated 
by $\{ i_{F(X)}(a) - (j_K \circ \phi)(a) ; a \in I_X\}$.  Then 
the Cuntz-Pimsner algebra ${\mathcal O}_X$ is the 
the quotient ${\mathcal T}_X/{\mathcal J}_X$ . 
Let $\pi : {\mathcal T}_X \rightarrow {\mathcal O}_X$ be the 
quotient map.  Put $S_{\xi} = \pi (T_{\xi})$ and 
$i(a) = \pi (i_{F(X)}(a))$. Let
$i_K : K(X) \rightarrow {\mathcal O}_X$ be the homomorphism 
defined by $i_K(\theta _{\xi,\eta}) = S_{\xi}S_{\eta}^*$. Then 
$\pi((j_K \circ \phi)(a)) = (i_K \circ \phi)(a)$ for $a \in I_X$.   
We note that  the Cuntz-Pimsner algebra ${\mathcal O}_X$ is 
the universal $C^*$-algebra generated by $i(a)$ with $a \in A$ and 
$S_{\xi}$ with $\xi \in X$  satisfying that 
$i(a)S_{\xi} = S_{\phi (a)\xi}$, $S_{\xi}i(a) = S_{\xi a}$, 
$S_{\xi}^*S_{\eta} = i((\xi | \eta)_A)$ for $a \in A$, 
$\xi, \eta \in X$ and $i(a) = (i_K \circ \phi)(a)$ for $a \in I_X$.
We usually identify $i(a)$ with $a$ in $A$.  We denote by 
${\mathcal O}_X^{alg}$ the $\ ^*$-algebra generated algebraically 
by $A$  and $S_{\xi}$ with $\xi \in X$. There exists an action 
$\gamma : {\mathbb R} \rightarrow Aut \ {\mathcal O}_X$
with $\gamma_t(S_{\xi}) = e^{it}S_{\xi}$, which is called the  
gauge action. Since we assume that $\phi: A \rightarrow L(X)$ is 
isometric, there is an embedding $\phi _n : L(X^{\otimes n})
 \rightarrow L(X^{\otimes n+1})$ with $\phi _n(T) = 
T \otimes id_X$ for $T \in L(X^{\otimes n})$ with the convention 
$\phi _0 = \phi : A \rightarrow L(X)$.  We denote by ${\mathcal F}_X$
the $C^*$-algebra generated by all $K(X^{\otimes n})$, $n \geq 0$ 
in the inductive limit algebra $\varinjlim L(X^{\otimes n})$. 
Let ${\mathcal F}_n$ be the $C^*$-subalgebra of ${\mathcal F}_X$ generated by 
$K(X^{\otimes k})$, $k = 0,1,\dots, n$, with the convention 
${\mathcal F}_0 = A = K(X^{\otimes 0})$.  Then  ${\mathcal F}_X = 
\varinjlim {\mathcal F}_n$.  
  
Let $p(z,w)$ be a  non-zero polynomial in two variables and 
 $\cC_p$ the algebraic correspondence of 
$p(z,w)$ on the Riemann sphere. 
Consider a $C^*$-algebra 
$A = C(\hat{\mathbb C})$ of 
continuous functions on $\hat{\mathbb C}$.   
Let $X = C(\cC_p)$. 
Then $X$ is a $A$-$A$ bimodule by 
$$
(a\cdot f \cdot b)(z,w) = a(z)f(z,w)b(w)
$$
for $a,b \in A$ and $f \in X$. We introduce an $A$-valued 
inner product $(\ |\ )_A$ on $X$ by 
\[
(f|g)_A(w)   = \sum_{\{\,z \in \hat{\mathbb C} \,|(z,w) \in \cC_p\,\}}
 e_p(z,w)\overline{f(z,w)}g(z,w)
\]
for $f,g \in X$ and $w \in \hat{\mathbb C}$.  We need branch index 
$e_p(z,w)$ in the formula of the inner product above.
Put $\|f\|_2 = \|(f|f)_A\|_{\infty}^{1/2}$.

\begin{lemma}
   The above $A$-valued inner product is well defined, 
that is, $\hat{\mathbb C} \ni w \mapsto (f|g)_A(w) 
\in {\mathbb C}$ is continuous.  
\label{inner-product}
\end{lemma}  

\begin{proof} 
If we consider 
$p(z,w) =   a_m(w)z^m + a_{m-1}(w)z^{m-1}+ \cdots + a_1(w)z + a_0(w)$, 
as a polynomial in $z$, then each coeficient $(a_k(w))_k$ is 
continuous in $w$. Then the continuity of 
the map $\hat{\mathbb C} \ni w \mapsto (f|g)_A(w) 
\in {\mathbb C}$ follows from the definition of the branch index and 
the continuity of the roots with multiplicities of a 
polynomial  on the  Riemann sphere. 
 See,  for example,  \cite{CC}.  

\end{proof} 
  
The left multiplication of $A$ on $X$ gives 
the left action $\phi : A \rightarrow L(X)$ 
such that 
$(\phi (a)f)(z,w) = a(z)f(z,w)$ for $a \in A$ and $f \in X$.

\begin{prop}
Let $p(z,w)$ be a  non-zero polynomial in two variables. 
Then  $X = C(\cC_p)$ is a full Hilbert $C^*$-bimodule over 
$A= C(\hat{\mathbb C})$ without completion. 
The left action $\phi : A \rightarrow L(X)$ is unital and 
faithful. 
\end{prop}

\begin{proof} Let $m$  be the degree of $p(z,w)$ in $z$. 
For any $f \in X = C(\cC_p)$, we have 
$$
\| f\|_{\infty} \leq \| f \|_2 
:= (\sup _w  \sum _{\{z \in  \hat{\mathbb C} |(z,w) \in \cC_p\}} 
e_p(z,w)|f(z,w)|^2)^{1/2}
\leq \sqrt{m} \| f\|_{\infty}. 
$$
Therefore two norms $\|\ \|_2$ and $\|\ \|_{\infty}$ are 
equivalent.  Since $C(\cC_p)$ is complete  with  respect to 
$\| \ \|_{\infty}$, it is also complete with  respect to 
$\|\ \|_2$.

Since $(1|1)_A(w) =  
\sum _{\{z \in \hat{\mathbb C} \ | \ (z,w) \in \cC_p\}} e_p(z.w)1 = m$, 
$(X|X)_A$ contains the identity $I_A$ of $A$.  Therefore 
$X$ is full. If $a\in A$ is not zero, then there exists 
$x_0 \in \hat{\mathbb C}$ with $a(x_0) \not = 0$.  Since the degree 
$m$ in $z$ of $p(z,w)$ is greater than or equal to one, there exists 
$w_0 \in \hat{\mathbb C}$ with $(x_0,w_0) \in \cC_p$. Choose 
$f \in X$ with $f(x_0,w_0) \not= 0$.  Then 
$\phi (a)f \not= 0$.  Thus $\phi$ is faithful.     

\end{proof}

\medskip\par
\noindent
{\bf Definition.} 
We introduce the $C^*$-algebra
${\mathcal O}_p(\hat{\mathbb C})$ associated with 
an algebraic correspondence 
$\cC_p = 
\{(z,w) \in \hat{\mathbb C} \times \hat{\mathbb C} \ | \ p(z,w) = 0 \}$
as a Cuntz-Pimsner algebra \cite{Pi} of the Hilbert $C^*$-bimodule 
$X_p= C(\cC_p)$ over 
$A = C(\hat{\mathbb C})$. 

\medskip\par 
A closed subset $J$ in $\hat{\mathbb C}$ is said to be  $p$-invariant 
if the following conditions are satisfied:
 For $z,w \in \hat{\mathbb C}$, 
 
(i)$z \in J$ and   
$p(z,w)=0$ implies $w \in J$, 

(ii)$w \in J$ and 
 $p(z,w)=0$ implies $z \in J$.  \\
\noindent 
Under the condition, we can define 
$\cC_p(J) = \{\,(z,w) \in J \times J \,\,|\,\,
p(z,w)=0\,\}$, $A= C(J)$, $X_p(J)= C(\cC_p(J))$ similarly. 
Then $X_p(J)$ is a full Hilbert $C^*$-bimodule ($C^*$-correspondence) over  
$A = C(J)$ and the left action is unital and faithful.  

We also introduce the $C^*$-algebra
${\mathcal O}_p(J)$ 
as a Cuntz-Pimsner algebra of the Hilbert $C^*$-bimodule 
$X_p(J)= C(\cC_p(J))$ over $A= C(J)$. 

\medskip\par
We define  the set $B(p)$ of "branched points" and 
the set $C(p)$ of "branched  values". 

\[
 B(p)  := \{\, z \in \hat{\mathbb C} \,\,|\,
\,  \text{there exists } w \in \hat{\mathbb C} 
\text{ such that }  p(z,w) = 0\ \text{ and }
e(z,w) \ge 2  \} . 
\]

\[
 C(p)  := \{\, w \in \hat{\mathbb C} \,\,|\,
\,  \text{there exists } z \in \hat{\mathbb C} 
\text{ such that }  p(z,w) = 0\ \text{ and }
e(z,w) \ge 2 \}.  
\]
In the above definitions, we may replace $e(z,w) \ge 2$ by 
$\frac{\partial p}{\partial z}(z,w) = 0$ after appropriate change 
of variables. 
Symmetrically we define 
\[
 \tilde{B}(p)  := \{\, w \in \hat{\mathbb C} \,\,|\,
\,  \text{there exists } z \in \hat{\mathbb C} 
\text{ such that }  p(z,w) = 0\ \text{ and }
\frac{\partial p}{\partial w}(z,w) = 0 \} .
\]

\[
 \tilde{C}(p)  := \{\, z \in \hat{\mathbb C} \,\,|\,
\,  \text{there exists } w \in \hat{\mathbb C} 
\text{ such that }  p(z,w) = 0\ \text{ and }
\frac{\partial p}{\partial w}(z,w) = 0 \}. 
\]

We need some known estimates of the above sets. 

\begin{lemma}
Let $p(z,w)$ be a  non-zero polynomial in two variables 
of degree $m$ in $z$ and degree $n$ in $w$. Then 
$B(p)$, $C(p)$, $\tilde{B}(p)$ and $\tilde{C}(p)$ 
are finte sets. More precisely  we have 
  $ ^{\#}B(p) \leq 2m(m-1)n$,  $ ^{\#}C(p) \leq 2(m-1)n$,  
$ ^{\#}\tilde{B}(p) \leq 2n(n-1)m$ and 
$ ^{\#}\tilde{C}(p) \leq 2(n-1)m$. 
\label{finiteset}
\end{lemma}  
\begin{proof}
It follows from proposition 2 in \cite{BP1} that 
$ ^{\#}C(p) \leq 2(m-1)n$.  Since $p(z,w)$ has degree 
$m$ in $z$, we also have $ ^{\#}B(p) \leq 2m(m-1)n$. 
The rest is symmetrically obtained.  
\end{proof}

Let  $I_X = I_{X_p(J)} = \phi^{-1}(\phi(C(J))\cap K(X_p(J)))$.

\begin{prop}
$I_{X_p(J)} 
= \{\, a \in C(J) \ | \  a|_{B(p)}\,\}=0$
\end{prop}  
\begin{proof}
This is a direct consequence of 
Proposition 4.4 in \cite{IKW} or \cite{MT}. 
\end{proof}

We consider Hilbert $C^*$-bimodules of iteration of the 
"algebraic  function". Put  
$X^{\otimes 2}_A =X \otimes_A X$, $X^{\otimes n}_A = X^{\otimes n-1}
\otimes_A X$.  

We define the path space  $\cP_n = \cP_n(J)$ of length $n$ in $J$ by 
$$
\cP_n = \{\,(z_1,z_2,\dots,z_{n+1}) \in J^{n+1}\,\,|\,\, 
(z_i,z_{i+1}) \in \cC_p(J)
,\,\,i=1,\dots,n \,\}.
$$
Then  $\cP_n$ is compact, since it is a continuous image of 
a compact subset. 
We extends branched index for paths of length $n$ as 
$$
e(z_1,z_2,\dots,z_{n+1}) =
e(z_1,z_2)e(z_2,z_3)\cdots e(z_{n},z_{n+1}). 
$$

Then ${\rm C}(\cP_n)$ is a Hilbert bimodule over $A$ by 
\begin{align*}
 & (a\cdot f \cdot b)(z_1,z_2,\dots,z_{n+1})  =
a(z_1)f(z_1,z_2,\dots,z_n)b(z_{n+1}) \\
 & (f|g)_A(w)  \\
& = \sum_{\{\,(z_1,z_2,\dots,z_n) | 
(z_1,z_2,\dots,z_n,w) \in
 \cP_n\,\}}e(z_1,z_2,\dots,z_n,w)
   \overline{f(z_1,z_2,\dots,z_n,w)}g(z_1,z_2,\dots,z_n,w)
\end{align*}
for $a$, $b \in A$, $f$, $g \in  C(\cP_n)$.  

\begin{lemma}
 The above $A$-valued inner product is well defined, 
that is, $\hat{\mathbb C} \ni z \mapsto (f|g)_A(w) 
\in {\mathbb C}$ is continuous, for any $f$, $g \in C(\cP_n)$. 
\label{lemma:continuous}
\end{lemma}  
\begin{proof}
It is enough to assume that $J = \hat{\mathbb C}$. 
We may and do assume that $n = 2$, because a similar argument 
holds for general $n$. We already know that $X \otimes_A X$ 
has an $A= C(\hat{\mathbb C})$-valued inner product. Therefore 
for $f_1\otimes f_2, g_1\otimes g_2 \in X \otimes_A X$, 
$\hat{\mathbb C} \ni w \mapsto 
(f_1\otimes f_2|g_1\otimes g_2)_A(w) \in {\mathbb C}$ is 
continuous. Put $f, g \in {\rm C}(\cP_2)$ by 
\[
f(z_1,z_2,w) = f_1(z_1,z_2)f_2(z_2,w), 
\ g(z_1,z_2,w) = g_1(z_1,z_2)g_2(z_2,w). 
\]
Then 
\begin{align*}
& (f|g)_A(w) = \sum_{\{\,(z_1,z_2)\in  \cP_1 | 
(z_1,z_2,w) \in
 \cP_2\,\}}e(z_1,z_2,w)
   \overline{f(z_1,z_2,w)}g(z_1,z_2,w) \\
 & = \sum_{\{\,(z_1,z_2)\in  \cP_1 | 
(z_1,z_2,w) \in
 \cP_2\,\}}e(z_1,z_2)e(z_2,w)
   \overline{f_1(z_1,z_2)f_2(z_2,w)}g_1(z_1,z_2)g_2(z_2,w) \\
& = \sum_{\{\,z_2 \in \hat{\mathbb C}  | 
p(z_2,w)=0 \,\}}e(z_2,w)
   \overline{f_2(z_2,w)}
(\sum_{\{\,z_1 \in  \hat{\mathbb C} | 
p(z_1,z_2)=0  \}}e(z_1,z_2)
   \overline{f_1(z_1,z_2)}g_1(z_1,z_2))g_2(z_2,w)) \\
& = \sum_{\{\,z_2 \in \hat{\mathbb C}  | 
p(z_2,w)=0 \,\}}e(z_2,w)
   \overline{f_2(z_2,w)}(f_1|g_1)_A(z_2)g_2(z_2,w) \\
& = (f_2|(f_1|g_1)_Ag_2)_A(w) 
= (f_1\otimes f_2|g_1\otimes g_2)_A(w). 
\end{align*}
Hence $w \mapsto (f|g)_A(w)$ is continuous. Then 
for  finte these sums 
\[
f(z_1,z_2,w) = \sum _if_{1,i}(z_1,z_2)f_{2,i}(z_2,w), 
\ g(z_1,z_2,w) = \sum _ig_{1,i}(z_1,z_2)g_{2,i}(z_2,w), 
\]
$w \mapsto (f|g)_A(w)$ is also continuous.
Put 
\[
C(\cP_2)^0 = \{f \in C(\cP_2) \ | 
 f(z_1,z_2,w) = \sum _{\text{finite} \ i} f_{1,i}(z_1,z_2)f_{2,i}(z_2,w)  
\text{ for } f_{1,i},f_{2,i} \in X \}.
\]

Since $C(\cP_2)^0$ is a $*-$subalgebra of $C(\cP_2)$ and 
separates two points, $C(\cP_2)^0$ is uniformly dense 
in $C(\cP_2)$. Note that uniformly norm 
$\| \centerdot \|_{\infty}$ and $\| \centerdot \|_2$ 
are equivalent, because 
$\| \centerdot \|_{\infty}  \leq \| \centerdot \|_2 
 \leq m^{n/2} \| \centerdot \|_{\infty}$. 
For any $f,g \in C(\cP_2)$, there exist 
sequences $(f_n)_n$ and $(g_n)_n$ in  $C(\cP_2)^0$ 
such that $f_n \rightarrow f$ and $g_n \rightarrow g$ uniformly. 
Since 
\begin{align*}
& |(f|g)_A(w) - (f_n|g_n)_A(w)|  \\
& \leq 
 \sum_{\{\,(z_1,z_2) | 
(z_1,z_2,w) \in
 \cP_2\,\}}e(z_1,z_2,w)
  | \overline{f(z_1,z_2,w)}g(z_1,z_2,w) - 
  \overline{f_n(z_1,z_2,w)}g_n(z_1,z_2,w)|, 
\end{align*}
we see that $(f_n|g_n)_A(w)$ converges to $(f|g)_A(w)$ 
uniformly in $w$.  Since a uniformly limit of continuous 
functions is continuous, $(f|g)_A(w)$ is continuous in $w$. 
\end{proof} 

Now it is clear to check the following propositon. 

\begin{prop}
Let $p(z,w)$ be a  non-zero polynomial in two variables. 
Then  $X = C(\cP_n)$ is a full Hilbert bimodule over 
$A= C(J)$ without completion. 
The left action $\phi : A \rightarrow L(X)$ is unital and 
faithful. 
\label{Pn}
\end{prop}

\begin{prop}
There exists an isometric $A$-$A$ bimodule homomorphism 
$\varphi : X^{\otimes n}_A \rightarrow {\rm C}(\cP_n)$ 
such that 
\[
 \varphi(f_1 \otimes f_2 \otimes \dots \otimes f_n)
(z_1,z_2,\dots,z_{n+1})
 = f_1(z_1,z_2)f_2(z_2,z_3)\cdots f_n(z_{n},z_{n+1})
\]
for $f_1, \cdots , f_n \in X$. 
\label{prop:continuous2}
\end{prop}  
\begin{proof}
It is easy to check $\varphi$ is a well defined 
 $A$-$A$ bimodule map. The proof of Lemma  
\ref{lemma:continuous} shows that $\varphi$ is isometric. 
Since $\| \centerdot \|_{\infty}$ and $\| \centerdot \|_2$ 
are equivalent, $\varphi$ is onto. 
\end{proof}

We need to define another  compact space  $\cG_n=\cG_n(J)$ by  
$$
\cG_n = \{\,(z_1,z_{n+1}) \in J^2\,\,
|\, \text{ there exists }(z_1,z_2,\dots,z_{n+1}) \in \cP_n\}.
$$
Then $C(\cG_n)$ is a Hilbert bimodule over $A$ by 

\begin{align*}
 & (a\cdot f \cdot b)(z_1,z_{n+1})  =
a(z_1)f(z_1,z_n)b(z_{n+1}) \\
 & (f|g)_A(w)   = \sum_{\{\,(z_1,z_2,\dots,z_n) | 
(z_1,z_2,\dots,z_n,w) \in
 \cP_n\,\}}e(z_1,z_2,\dots,z_n,w)
   \overline{f(z_1,w)}g(z_1,w)
\end{align*}
for $a$, $b \in A$, $f$, $g \in  C(\cG_n)$.  Define 
a continuous onto map $\rho : \cP_n \rightarrow \cG_n$ 
by $\rho((z_1,z_2,\dots,z_{n+1})) = (z_1,z_{n+1})$ for 
$(z_1,z_2,\dots,z_{n+1}) \in \cP_n$. The it is clear that 
the induced map  
$\rho^* : C(\cG_n) \rightarrow C(\cP_n) $ 
defined by $\rho^*(f) = f \circ \rho$ is an 
isometric Hilbert bimodule embedding.  

\section{simplicity and pure infiniteness}

In this section we consider a sufficient condition for 
a polynomial so that the associated ${\mathcal O}_p(J)$ 
is simple and purely infinite. 

Let $J$ be a  $p$-invariant subset of  $\hat{\mathbb C}$. 
For any subset $U$ of $J$ and a natural number $n$,  
we define a subset $U^{(n)}$ of $J$ by 
\[
 U^{(n)} = \{\,w \in J \,\,| \,\,(z_1,z_2,\dots,z_n,w)\in \cP_n
 \quad \text{for some} \quad z_1 \in U,\, z_2 \dots z_n \in J \,\}.
\]
\medskip\par
\noindent
{\bf Definition.} 
Let $p(z,w)$ be a  non-zero polynomial in two variables and 
$J$ a  $p$-invariant subset of  $\hat{\mathbb C}$. 
Then $p$ is said to be {\it expansive} on $J$ if 
for any nonempty open set $U \subset J$ in $J$ with the relative 
topology 
there exists a natural number $n$ 
such that   $U^{(n)} = J$. 

\medskip\par
\noindent
{\bf Example.} 
Let $R(z) = \frac{P(z)}{Q(z)}$ be a rational function with 
polynomials $P(z)$, $Q(z)$ and $\deg R \geq 2$. 
Put $p(z,w)=Q(z)w-P(z)$. Then  $U^{(n)}$ is exactly $R^n(U)$. 
Therefore 
$p$ is expansive 
on the Julia set $J_R$ by \cite[Theorem 4.2.5]{Be}. 

\medskip\par
\noindent
{\bf Example.} Let $p(z,w) = z^2 + w^2 -1$. Then 
$J := \{ 0, 1, -1 \}$ is a $p$-invariant set and 
$p$ is not expansive on $J$. In fact, let $U = \{0\}$, then 
$U^{(2n)} = \{0\}$ and $U^{(2n+1)} = \{1,-1\}$. 
And $p$ is not expansive  
on $\hat{\mathbb C}$, because an open set 
$U := \hat{\mathbb C} \setminus \{ 0, 1, -1 \}$ is 
$p$-invariant and $U^{(n)} = U \not = \hat{\mathbb C}$ 
for any $n$.  In general, for any polynomial, 
if $p$ has a finite $p$-invariant 
set, then $p$ is not expansive on $\hat{\mathbb C}$ similarly. 

\medskip\par
\noindent
{\bf Example.} Let $p(z,w) = z^m - w $, $m \geq 2$. 
Then $J := {\mathbb T}$ is a $p$-invariant set and 
$p$ is expansive on $J$, because ${\mathbb T}$ is a 
Julia set of $w = R(z) = z^m$. 

Let $p(z,w) = z- w^n $, $n \geq 2$. 
Then $J := {\mathbb T}$ is a $p$-invariant set but 
$p$ is not expansive on $J$. In fact 
$U := \hat{\mathbb C} \setminus \{1 \}$ 
is $p$-invariant and $U^{(k)} = U \not = \hat{\mathbb C}$ 
for any $k$.

More generaly we have a criterion. 

\begin{prop}
Let $p(z,w) = z^m - w^n$ for natural 
numbers $m$ and $n$. 
Then $J := {\mathbb T}$ is a $p$-invariant set. 
And  
$p$ is expansive on ${\mathbb T}$ 
if and only if $n$ is not devided by $m$. 
\label{expansive}
\end{prop}
\begin{proof}
Suppose that $n$ is devided by $m$, so that $n = mj$ for some 
$j \in {\mathbb N}$. 
Let $U = \{z \in {\mathbb T} 
\ | \ z^m \not= 1 \}$. Then for any $k \in {\mathbb N}$, 
$1$ is not in $U^{(k)}$. In fact, if $1$ were in $U^{(k)}$, 
then there exist $(z_1,z_2,\dots,z_k,1)\in \cP_k $ such that 
$z_1 \in U$. Hence $z_k^m = 1$ and 
$z_{k-1}^m = z_k^n = z_k^{mj} = 1$.  We continue this argument 
to obtain $z_1^m = 1$. This contradicts the fact that $z_1 \in U$. 
Therefore  $p$ is not expansive on ${\mathbb T}$.  

Nextly, suppose that $n$ is not devided by $m$. Let $d$ be the 
greatest common divisor of $m$ and $n$. Then $m = m_0d$ and 
$n = n_0d$ for some natural numbers $m_0$ and $n_0$. Since 
$n$ is not devided by $m$, $m_0$ is greater than or equal to 2. 
We identify ${\mathbb T}$ with ${\mathbb R} \  (\mod 1)$ by 
$z = e^{2\pi i\alpha}$ and $w = e^{2\pi i\beta}$. Then 
$z^m - w^n = 0$ means that $m\alpha = n\beta - k$ for some 
integer $k$. Hence 
\[
\cC_p(J) \cong \{\,([\alpha],[\beta]) 
\in {\mathbb R}/{\mathbb Z}\times {\mathbb R}/{\mathbb Z} \,\,|\,\,
\beta = \frac{m}{n} \alpha + \frac{k}{n} \text{ for some integer } k \,\}.
\]
Then $\cC_p(J)$ has $d$ connected components, because 
$m{\mathbb R} = d{\mathbb R} \ (\mod n)$. 

For $([\alpha],[\beta]) \in \cG_2(J)$, there exist $k_1,k_2 \in {\mathbb Z}$ 
such that 
\[
\beta = \frac{m}{n}( \frac{m}{n} \alpha + \frac{k_1}{n}) + \frac{k_2}{n}
= \frac{m^2}{n^2} \alpha + \frac{mk_1 + nk_2}{n^2}. 
\]
Since $m \mathbb Z + n \mathbb Z = d\mathbb Z$, 
$([\alpha],[\beta]) \in \cG_2(J)$ if and only if 
there exists $k \in {\mathbb Z}$ 
such that 
$\beta = \frac{m^2}{n^2} \alpha + \frac{dk}{n^2}$
We continue in this way to get that 
$([\alpha],[\beta]) \in \cG_r(J)$ if and only if 
there exists $k \in {\mathbb Z}$ 
such that 
$\beta = \frac{m^r}{n^r} \alpha + \frac{d^{r-1}k}{n^r}$. 
Since  $m^r {\mathbb Z} + n^r \mathbb Z = d^r{\mathbb Z}$, 
$\cG_r(J)$ has $\frac{d^n}{d^{n-1}}= d$ connected components.
To avoid overlapping, 
we see only one connected component. Hence 
we need to cover an interval 
$[0, \frac{d^{r-1}}{n^r}d] = [0, \frac{d^r}{n^r}]$. 
Take an open interval $I = (0,\frac{1}{m_0^r} + \varepsilon )$. Since
$\frac{m^r}{n^r}\frac{1}{m_0} = \frac{d^r}{n^r}$, 
$I^{(r)}$ contains 
\[
(0,\frac{d^r}{n^r}] \cup (\frac{d^r}{n^r}, \frac{2d^r}{n^r}] \cup \dots 
\cup (\frac{(n_0^r -1)d^r}{n^r}, \frac{n_0^rd^r}{n^r}]
\cup \{0\} = [0,1] \  (\mod {\mathbb Z}).
\]
It is also true if we replace $I$ by a translation of $I$. 
Now for any open set $U \subset J$, there exists an open interval 
$(a,b)$ with $(a,b) \subset U$. Choose a natural number $r$ 
such that $|b-a| > 1/{m_0^r}$.  Then by the preceeding argument 
we see that $(a,b)^{(r)} = [0,1]$. 
Hence $U^{(r)} = [0,1] (mod \ {\mathbb Z}) = J$.  This shows that 
$p$ is expansive on ${\mathbb T}$.  

\end{proof}

\medskip\par
\noindent
{\bf Definition.}
Let $N$ be a natural number. We define the set ${\rm GP}(N)$
of $N$-{\it generalized periodic points}  by 
\begin{align*}
     {\rm GP}(N)
=  & \{ w \in J \,|\, \exists z \in J  \ \exists m,\,n  \quad0 \le m \ne  n \le N,
  \exists (z,z_2,z_3,\dots,z_n,w) \in \cP_n,\, \\
 &  \exists (z,u_2,u_3,\dots,u_m,w) \in \cP_m \, \}. 
\end{align*}
Let  $R(z) = \frac{P(z)}{Q(z)}$ be a rational function with 
polynomials $P(z)$, $Q(z)$. Put $p(z,w)=Q(z)w-P(z)$. Then 
\[
{\rm GP}(N) = \cup _{n=1}^N 
\{ w \in \hat{\mathbb C} \ | \ R^n(w) = w \}.  
\]
In fact, if $R^n(w) = w$ for some $n \leq N$, then 
it is clear that $w \in  {\rm GP}(N)$. Conversely 
let $w \in {\rm GP}(N)$, then there exists 
$z$ such that $w = R^n(z) = R^m(z)$ for 
some $0 \leq m < n \leq N$. Then $R^{n-m}(w) = w$. 

\medskip\par
\noindent
{\bf Definition.}
A polynomial $p$ in two variables
is said to be {\it free} on $J$ if for any natural number $N$,  
${\rm GP}(N)$ is a finite set.  

For example, let  $R(z) = \frac{P(z)}{Q(z)}$ be a rational function
with polynomials $P(z)$, $Q(z)$. Put $p(z,w)=Q(z)w-P(z)$. If 
$\deg R \geq 2$, then $p$ is free on any $p$-invariant set $J$.

\begin{lemma} Let $p(z,w)=z^m-w^n$.  Then 
$p$ is free on  $J = {\mathbb T}$ 
if and only if   $m \ne n$,.  
\end{lemma}
\begin{proof}Assume that $m \ne n$. 
We identify ${\mathbb T}$ with ${\mathbb R}/{\mathbb Z}$ by 
$z = e^{2\pi i\alpha}$ and $w = e^{2\pi i\beta}$. 
For any natural number $N$,  
$[\beta] \in {\rm GP}(N)$ if and only if there exist
$[\alpha] \in {\mathbb R}/{\mathbb Z}$ and 
$ r,\,s \quad0 \le r \ne  s \le N$ such that 
$([\alpha],[\beta]) \in \cG_r(J)$ and 
$([\alpha],[\beta]) \in \cG_s(J)$. Therefore 
there exist $k_1, k_2 \in {\mathbb Z}$ with 
$0 \leq k_1 \leq n^r-1$ and $0 \leq k_2 \leq n^s-1$ 
such that 
\[
\beta = \frac{m^r}{n^r} \alpha + \frac{d^{r-1}k_1}{n^r}
      = \frac{m^s}{n^s} \alpha + \frac{d^{s-1}k_2}{n^s}. 
\]
Since two segment with different slopes meet at most one 
point, $ {}^{\#}{\rm GP}(N) \leq n^{3N}$.  
Hence $p$ is free on ${\mathbb T}$. 

Conversely assume that $m = n$. Then any 
$(z,z,\dots,z) \in J^{k+1}$ is in  $\cP_k(J)$. Hence 
for any natural number $N$,  
${\rm GP}(N) = {\mathbb T}$ is an infinite set. Thus 
$p$ is not free on ${\mathbb T}$. 
\end{proof}

\medskip\par
\noindent
{\bf Remark.} The above example is related with an 
example by Katsura in \cite[section 4]{Ka2} .
 If $m$ and $n$ are relatively 
prime, then his example coincides with our example. If
 $m$ and $n$ are not relatively 
prime, then his example is different with ours. 
In fact our $\cP_n({\mathbb T})$ is not 
connected if $m$ and $n$ are not relatively 
prime.  But they are isomorphic as bimodule.

\begin{prop} Let $R_i(z)=P_i(z)/Q_i(z)$ $i=1,\dots, r$ 
be rational functions with polynomials $P_i(z)$, $Q_i(z)$. 
Put $p(z,w) = (Q_1(z)w-P_1(z))\cdots (Q_r(z)w-P_r(z))$. 
Let $J \subset \hat{\mathbb C}$ be a $p$-invariant 
closed subset. Assume that each $\deg R_i \geq 2$ and 
$\deg R_1, \dots , \deg R_r$ are relatively prime, then 
$p$ is free on $J$.  Furthermore, if 
$J$ is a Julia set for some $R_i$, then 
$p$ is expansive on $J$. 
\end{prop} 
\begin{proof} Let  $N$  be a natural number and $m, n$ integers with 
$0 \leq n < m \leq N$.  
For $i_1, \dots , i_m, j_1, \dots , j_n = 1,2, \dots, r$, 
we shall show that 
\[
\  M := ^{\#}\{ z \in \hat{\mathbb C} \ 
| \ R_{i_m}\circ \dots \circ R_{i_1}(z) = 
    R_{j_n}\circ \dots \circ R_{j_1}(z) \} < \infty . 
\]
On the contrary, assume that $M = \infty$. 
Then covering degrees of both sides coincide.  Count the 
covering derees and rearrange them. Then we  have 
\[
(\deg R_1)^{s_1} \dots (\deg R_r)^{s_r} =  
(\deg R_1)^{t_1} \dots (\deg R_r)^{t_r} 
\]
with $s_1 + \dots + s_r = m$ and 
$t_1 + \dots + t_r = n$. Since 
$\deg R_1, \dots , \deg R_r$ are relatively prime, 
$s_i = t_i$ for $i = 1, \dots , r$.  Then $m = n$.  
This contradicts the fact that $n < m$.  Hence  $M <\infty$. 
Therefore $Q(m,n) := \{ z \in \hat{\mathbb C} \ | \ 
\text{ there exists } w \in \hat{\mathbb C} \text{ such that } 
(z,w) \in \cG_m, \ (z,w) \in \cG_n \}$ is a finite set. Hence   
\[ 
{\rm GP}(N)
=  \{ w \in J \,|\, \exists z \in J  \ \exists m,\,n  \quad0 \le n < m \le N,
  \exists (z,w) \in \cG_m,\, 
  \exists (z,w) \in \cG_n \, \}
\]
is also a finite set. 
This shows that  $p$ is free on $J$. 

It is evident that, if $J$ is a Julia set for some $R_i$, then 
$p$ is expansive on $J$.

\end{proof}

\medskip\par
\noindent
{\bf Example.} 
Let  $m$ and $n$ be natural numbers and relatively prime. 
Consider  $p(z,w)=(w-z^m)(w-z^n)$. We note that 
 $J={\mathbb T}$ is the common Julia set of
 $w = z^m$ and $w = z^n$. 
Then $p$ is free on $J$ and expansive on $J$. 
We note that there appears a new branched point $(1,1)$ in $\cC_p$. 

\medskip\par
\noindent
{\bf Example.} Let $R_1(z)=\frac{(z^2 + 1)^2}{4z(z^2 - 1)}$ 
be a  rational function given by Lattes. Then the Julia set
$J_{R_1} = \hat{\mathbb C}$. Let 
 $R_2(z)=P_2(z)/Q_2(z)$ be any rational function 
with odd degree. 
Put $p(z,w) =
 ((4z(z^2 - 1))w-(z^2 + 1)^2)(Q_2(z)w-P_2(z))$. 
Let $J=\hat{\mathbb C}$. Then  $p$ is expansive on $J$ 
and free on $J$. 

\medskip\par
\noindent
{\bf Example.} 
Let $i_1$, $\cdots$, $i_n$, $j_1$, $\cdots$, $j_n$ be natural 
numbers. Assume that $i_k \not= 1$ or $j_k \not= 1$ 
for each $k$.  Suppose that those which are not equal to $1$ 
are relatively prime. Put $J={\mathbb  T}$. Let  
\[
 p(z,w) = (z^{i_1}-w^{j_1})(z^{i_2}-w^{j_2})\cdots (z^{i_n}-w^{j_n}). 
\]
Then $p$ is free on $J$.  

\medskip\par
\noindent
{\bf Example.} 
Let $m$ be a natural number with $m \geq 2$. Put
$p(z,w) = (w-z^m)(w^m-z)$. Then $p$ is not 
free on ${\mathbb T}$. In fact, there exist different paths 
$(z,z^m,z,z^m,z) \in \cP_4(\mathbb T)$ and 
$(z,z^m,z) \in \cP_2(\mathbb T)$. Hence 
${\rm GP}(4) = \mathbb T$. 

\medskip\par
\noindent
{\bf Example.} Let $p(z,w) = z^2 + w^2 -1$.  Then 
$p$ is not free on $J = {\hat{\mathbb C}}$. 
In fact, choose any 
$(z,w) \in \cC_p$. Then there exist different paths 
$(z,w,z,w,z) \in \cP_4(\hat{\mathbb C})$ and 
$(z,w,z) \in \cP_2(\hat{\mathbb C})$. Hence 
${\rm GP}(4) = \hat{\mathbb C}$. 


\begin{lemma}Suppose that $p$ is expansive on a $p$-invariant 
subset $J$. Then 
for any non-zero positive element $a \in A$ and for any  
$\varepsilon > 0$ there exist $n \in \mathbb{N}$ and 
$f \in X^{\otimes n}$  with $(f|f)_A = I$ such that  
$$
 \|a\| -\varepsilon \le S_f^* a S_f \le \|a\|. 
$$
\label{lemma:epsilon}
\end{lemma}
\begin{proof} Let $x_0$ be a point in $J$ with 
$|a(x_0)| = \| a \|$.  For any  
$\varepsilon > 0$ there exist an open neighbourhood $U$ of $x_0$
in $J$ such that for any $x \in U$  we have 
$\|a\| -\varepsilon  \le a(x) \le \|a\| $. 
Choose  another open neighbourhood $V$ of $x_0$ in $J$ and a compact 
subset $K \subset J$ satisfying $V \subset K \subset U$.
Since $p$ is expansive on $J$, there exists $n \in \mathbb{N}$ such that 
$V^{(n)} =J$. 
We identify $X^{\otimes n}$ with $C(\cP_n) \supset \rho^*(C(\cG_n))$
as in the paragraph after Proposition \ref{prop:continuous2}. 
Define closed subsets $F_1$ and $F_2$ of $J\times J$ by  

\begin{align*}
 F_1 & = \{(z,w)\in J\times J | (z,w) \in \cG_n, z \in K\},  \\
 F_2 & = \{(z,w) \in J\times J| (z,w) \in \cG_n, z \in U^{c} \}. 
\end{align*}

Since $F_1 \cap F_2 = \phi $, there exists  
$g \in C(\cG_n)$ such that 
$0 \le g(z,w) \le 1$  and 
$$
 g(z,w) =\left\{\begin{array}{cc}
  1,  & (z,w)\in F_1  \\
  0,  & (z,w)\in F_2 . 
       \end{array}\right.
$$
\par
Since $V^{(n)}=J$, for any $w\in J$ there exists 
$z_1 \in V$ such that $(z_1,w) \in \cG_n$. Then 
$(z_1,w) \in F_1$ and $g(z_1,w) = 1$. Therefore 
 
\begin{align*}
 (\rho^*(g)|\rho^*(g))_A(w) 
& = \sum_{\{(z_1,\dots,z_n) \in \cP_{n-1}|(z_1,\dots,z_n,w) \in \cP_{n}  \} }
e(z_1,\dots,z_n,w) |g(z_1,w)|^2 \\
            & \ge |g(z_1,w)|^2  =1 . 
\end{align*}

Let $b:=(\rho^*(g)|\rho^*(g))_A$.  
Then $b(y) = (\rho^*(g)|\rho^*(g))_A(y) \ge 1$.  Thus 
$b \in A$ is positive and invertible.
We put $f := \rho^*(g) b^{-1/2} \in X^{\otimes n}$.  Then  
$$
 (f|f)_A  = b^{-1/2} (g|g)_A b^{-1/2} 
           = I.
$$

For any $w \in J$ and $(z_1,w) \in \cG_n$, if $z \in U$, 
then $\|a\|-\varepsilon \le a(z)$, and  
if $z \in U^c$, then  $f(z_1,\dots,z_n,w)=g(x,y) b^{-1/2}(w) = 0$. 
Therefore   
\begin{align*}
 \|a \|-\varepsilon 
& = (\|a\|-\varepsilon) (f|f)_A(y) \\
& = (\|a\|-\varepsilon) 
\sum_{\{(z_1,\dots,z_n) \in \cP_{n-1}|(z_1,\dots,z_n,w) \in \cP_{n}  \}} 
 e(z_1,\dots,z_n,w) |f(z_1,\dots,z_n,w)|^2             \\
& \leq 
\sum_{\{(z_1,\dots,z_n) \in \cP_{n-1}|(z_1,\dots,z_n,w) \in \cP_{n}  \}} 
 e(z_1,\dots,z_n,w) a(z_1)|f(z_1,\dots,z_n,w)|^2 \\
& = (f|af)_A(w) = S_f^* a S_f (w). 
\end{align*}
It is clear that
$S_f^* a S_f = (f|af)_A \leq \| a \| (f|f)_A = \| a \| $.
\end{proof}

\begin{lemma}
Suppose that $p$ is expansive on a $p$-invariant 
subset $J$. Then 
for any non-zero positive element $a \in A$ and for any  
$\varepsilon > 0$ with $0 < \varepsilon < \|a\|$, 
there exist $n \in \mathbb{N}$ and $u \in X^{\otimes n}$ 
such that  
$$
 \| u \| _2 \le (\|a\| - \varepsilon)^{-1/2} \qquad \text{and} 
\quad  S_u^* a S_u = I. 
\label{lemma:u-epsilon}
$$
\end{lemma}
\begin{proof} 
This is exactly as same as \cite[Lemma 3.5]{KW1}.
\end{proof}

A step in the proof of the main theorem is to 
show a certain element $S$ in a Cuntz-Pimsner algebra is 0. It is 
enough to show a corresponding element $T$ 
in the Toeplitz algebra is 0.  Since the Toeplitz algebra acts 
on the Fock module and the Fock module is realized as a function 
space, we can calculate $Tx = 0$ concretely.  

We write $A=X^{\otimes 0}$. If $a \in A$, then 
$T_a$ means $\phi(a)\otimes I_n$ on $X^{\otimes n}$. 
The following lemma is a key of the proof of our main theorem. 

\begin{lemma} Let $i$ and $j$ be integers 
with  $i, j  \geq 0$ and $i \ne j$. 
Take $x \in X^{\otimes i}$ and $y \in X^{\otimes j}$. 
Suppose that $a \in A = C(J)$ satisfies the following 
condition:

 For any $(z_1,z_2,\dots,z_i,w) \in \cP_i$, 
$(u_1,u_2,\dots,u_j,w) \in \cP_j$, we have $a(z_1)a(u_1)=0$. \\
Then we have $aT_xT_y^*a^*=0$.
\label{lem:Fock}
\end{lemma}
\begin{proof} It is enough to show $T_{ax}T_{ay}^*f=0$ 
for any $f \in X^{\otimes r}$, $r = 0,1,2,\dots  $. 
If $r < j$, then $T_{ax}T_{ay}^*f=T_{ax}0=0$. Hnece 
we may assume that  $r \geq j$ and $ f = f_1 \otimes f_2$ for 
 $f_1 \in X^{\otimes j}$,  $f_2 \in X^{\otimes (r-j)}$. 
\begin{align*}
 &  (T_{ax}T_{ay}^*)(f_1 \otimes
 f_2)(z_1,z_2,\dots,z_{i},z_{i+1},\dots, z_{i +r-j + 1})
 \\
= & (T_{ax} (ay|f_1)_Af_2)(z_1,z_2,\dots,z_{i},z_{i+1},\dots, z_{i +r-j + 1}) \\= & (ax \otimes (ay|f_1)_Af_2)(z_1,z_2,\dots,z_{i},z_{i+1},\dots, z_{i +r-j + 1}) \\
= & a(z_1)x(z_1,\cdots,z_{i+1})
 (ay|f_1)_A(z_{i+1})f_2(z_{i+1},\dots,z_{i +r-j+1}) \\
= & a(z_1)x(z_1,\cdots,z_{i+1})  \\
 & \cdot \left(  \sum_{(u_1,\dots,u_j,z_{i+1}) \in \cP_j}
e(u_1,\dots,u_j,z_{i+1})
   \overline{a(u_1)y(u_1,\cdots,u_j,z_{i+1})}
f_1(u_1,\cdots,u_j,z_{i+1})\right)
 \cdot  \\
 &    f_2(z_{i+1},\dots,z_{i +r-j+1}) \\
= & a(z_1)\overline{a(u_1)}x(z_1,\cdots,z_{i+1})  \\
 & \cdot \left(  \sum_{(u_1,\dots,u_j,z_{i+1}) \in \cP_j}
e(u_1,\dots,u_j,z_{i+1})
   \overline{y(u_1,\cdots,u_j,z_{i+1})}
f_1(u_1,\cdots,u_j,z_{i+1})\right)
 \cdot  \\
 &    f_2(z_{i+1},\dots,z_{i +r-j+1}) \\
 =& 0. 
\end{align*}

\end{proof}

We need to prepare the following elementary fact:

\begin{lemma}
Suppose that $p(z_0,w_0) = 0$ , 
$\frac{\partial p}{\partial z}(z_0,w_0) \not= 0$ and 
$\frac{\partial p}{\partial w}(z_0,w_0) \not= 0$. 
Then there exist an open neibourhood $U$ of $z_0$, 
an open neibourhood $V$ of $w_0$ and homeomorphism 
$\varphi : U \rightarrow V$ such that for any 
$z \in U$, $w \in V$, $p(z,w) = 0$ if and only if 
$w = \varphi (z)$.  
\label{lemma:implicit}
\end{lemma}

\begin{lemma}Assume that 
$p$ is free on $J$.  Suppose that $J$ has no isolated points. 
Let $N$ be a natural number. Then for any non-empty open set $U$ in $J$, 
there exist points $w_0 \in U, z_i \in J$ \ ($i = 1, \dots, m^N$), 
an open neibourhood $V$ of $w_0$ with  $V \subset U$, 
open neibourhoods $W_i$ of $z_i$ in $J$ and  
homeomorphisms $\Phi_i: W_i \rightarrow V$ for 
$i = 1, \dots, m^N$ satisfying the following. 
\begin{enumerate}
 \item $W_i \cap W_j = \empty$ for $i \not= j$.  
 \item For any  $z \in W_i $, $w \in V$, we have 
     $(z,w) \in \cG_N$ 
        if and only if $w = \Phi_i (z)$, in particular 
$w_0 = \Phi_i (z_i)$
 \item  For any $s\in J$ with $(z_i,s) \in \cG_k$ 
      for some $k$  ($1 \leq k \leq N)$, there exist an 
open neibourhood $W_{i,s}$ of $s$ 
and homeomorphisms $\Phi_{i,s}: W_i \rightarrow W_{i,s}$ satisfying 
the following : for any  $z \in W_i $, $w \in W_{i,s}$, we have 
     $(z,w) \in \cG_k$ 
        if and only if $w = \Phi_{i,s} (z)$. 
 \item These open neibourhoods  $W_i$ and $W_{i,s}$\ for $i,s$ have 
       empty intersection each other. 
\end{enumerate} 
\label{lemma:free}
\end{lemma}
\begin{proof}Let 
$D_1$ be the set of $w \in J$ satisfying that
there exist $u \in J, \ z \in {\rm GP}(N)$  
sucn that $(u,w) \in \cG_N, \ (u,z) \in  \cG_k$  \ 
for some $k = 0,1,\dots, N$.  

Since $p$ is free on $J$, ${\rm GP}(N) $ is a finte set. Hence 
$D_1$ is also a finte set. 
Consider the set $D_2$ of $ w \in J$ satisfying that
there exist $u,z \in J$  
such that  $(u,w) \in \cG_N, \ (u,z) \in  \cG_k $
for some $k = 0,1,\dots, N $ and 
$z$  is in $B(p)$, $C(p)$, $\tilde{B}(p)$ or  $\tilde{C}(p)$.  
Then $D_2$ is a finite set. Since $D_1 \cup D_2$ is a finite set and 
$J$ has no isolated points, there exist a non-empty open set 
$V_0 \subset U$  such that 
$V_0 \subset U \setminus (D_1 \cup D_2)$. Choose 
$w_0 \in V_0 \subset U \setminus (D_1 \cup D_2)$.
There exist distinct $z_i \in J $  for  $i = 1, \dots, m^N$
such that $(z_i,w_0) \in \cG_N$. By 
the Lemma \ref{lemma:implicit} , we can choose 
a sufficiently small non-empty open set  $V \subset V_0$, non-empty 
open neibourhoods $W_i$ of $z_i$ and 
homeomorphisms $\Phi_i: W_i \rightarrow V$ for 
$i = 1, \dots, m^N$ satisfying all the above requirements. 
\end{proof}

\begin{prop}Let $J$ be a $p$-invariant set  with no isolated points.  
Suppose that $p$ is expansive and  free on $J$. 
For $N \in \mathbb{N}$, for any $T \in L(X^{\otimes N})$, 
for any $\varepsilon > 0$,  there exists 
$a \in A^+  =C(J)^+ $ with $\|a\|=1$  such that 
\begin{align*}
 \|\phi(a) T\|^2  & \ge \|T\|^2 - \varepsilon,  \\
 a S_xS_y^*a  & = 0 \quad \text{ for any } x \in X^{\otimes i}, 
   \text{ for any }  y \in
 X^{\otimes j}, 0 \le i,j \le N, i \ne j. 
\end{align*}
\label{prop:free}
\end{prop}
\begin{proof}
For $N \in \mathbb{N}$, for any $T \in L(X^{\otimes N})$, 
for any $\varepsilon > 0$,  there exists 
$f \in X^{\otimes N}$ with $\|f\|_2 = 1$ such that 
$\|T\|^2 \ge \|Tf \|_2^2 > \|T\|^2-\varepsilon$. 
Hence there exists $w_1 \in J$ such that
\[
 \|Tf\|_2^2 = \sum_{\{\,(z_1,z_2,\cdots,z_N) | (z_1,z_2,\cdots,z_N,w_1) \in
 \cP_N\,\}} e(z_1,z_2,\cdots,z_N,w_1)|(Tf)(z_1,z_2,\cdots,z_N,w_1)|^2
\]
Since the function 
\[
 w \mapsto \sum_{\{\,(z_1,z_2,\cdots,z_N) | (z_1,z_2,\cdots,z_N,w) \in
 \cP_N\,\}} e(z_1,z_2,\cdots,z_N,w) |(Tf)(z_1,\dots,w)|^2
\]
is continuous, there exists an open neibourhood $U$ of $w_1$ such that 
for any $w \in U $ 
\[
  \sum_{\{\,(z_1,z_2,\cdots,z_N) | (z_1,z_2,\cdots,z_N,w) \in
 \cP_N\,\}}e(z_1,z_2,\cdots,z_N,w) |(Tf)(z_1,\cdots,z_N,w)|^2 > \|T\|^2
 - \varepsilon
\]
By  Propostion \ref{prop:free}, 
there exist points $w_0 \in U, z_i \in J$ \ ($i = 1, \dots, m^N$), 
an open neibourhood $V$ of $w_0$ with  $V \subset U$, 
open neibourhoods $W_i$ of $z_i$ in $J$ and  
homeomorphisms $\Phi_i: W_i \rightarrow V$ for 
$i = 1, \dots, m^N$ satisfying the conditions in the lemma. 
Choose $b \in A = C(J)$ satisfying 
\[
 b(w_0)=1, \quad 0 \le b(w) \le 1, \quad \supp\, b  \subset V. 
\]
Define $a \in C(J)$ by 
\[
        a(z) =
       \begin{cases}
	b(\Phi_i (z)) \qquad z \in W_i \\
        0       \qquad {\rm otherwize}. 
       \end{cases}
\]
Then this function $a$ satisfies the condition in Lemma \ref{lem:Fock}. 
Therefore 
for any $x \in X^{\otimes i}$, $y
\in X^{\otimes j}$, $0 \le i,j \le N$, $i \ne j$ 
we have $aS_xS_y^*a^*=0$.  

Moreove we have  
\begin{align*}
 \|\phi(a)Tf\|_2^2 & = \sup_w \sum_{\{\,(z,z_2,\dots,z_N) |
 (z,z_2,\dots,z_N,w) \in \cP_N\,\}} e(z,z_2,\dots,z_N,w)
 |a(z)(Tf)(z,\dots,w)|^2 \\
        & \ge  \sup_w \sum_{\{\,(z,z_2,\dots,z_N) | (z,z_2,\dots,z_N,w) \in
 \cP_N\,\}}
 e(z,z_2,\dots,z_N,w)|(Tf)(z,\dots,w) b(w)|^2 \\
                   & \ge  \sum_{\{\,(z,z_2,\dots,z_N) | (z,z_2,\dots,z_N,w_0)
        \in \cP_N\,\}} e(z,z_2,\dots,z_N,w_0)|(Tf)(z,\dots,w_0)b(w_0)|^2 \\
                   & >  \|T\|^2 - \varepsilon. 
\end{align*}
\end{proof}

It is important to recall the fact that 
there exists an isometric $\ ^*$-homomorhism
\[
\varphi:  L(X^{\otimes N}) \supset A \otimes I^N + K(X) \otimes I^{N-1} +
\cdots + K(X^{\otimes N}) \rightarrow \cO_{p}(J)^{\bf T}
\]
as in Pimsner 
\cite[Proposition 3.11]{Pi}  and  
Fowler-Muhly-Raeburn \cite[Proposition 4.6]{FMR}
such that 
\[
\varphi ( a + \theta_{x_1 \otimes \dots \otimes x_k, 
         y_1 \otimes \dots \otimes y_k}) 
    = a +  S_{x_1} \dots S_{x_k}S_{y_k}^* \dots S_{y_1}^*.
\]
To simplify notation, we put $S_x = S_{x_1} \dots S_{x_k}$ 
for $x = x_1 \otimes \dots \otimes x_k \in X^{\otimes k}$ .

\begin{lemma}Let $J$ be a $p$-invariant set  with no isolated points.  
Suppose that $p$ is expansive and  free on $J$. 
Let $b = c^*c$ for some $c \in {\mathcal O}_X^{alg}$.  
We decompose  $b = \sum _j b_j$ with 
$\gamma _t(b_j) = e^{ijt}b_j$.
For any  $\varepsilon >0 $
there exists $P \in A$ with $0\le P \le I$ satisfying the 
following:  
\begin{enumerate}
 \item $Pb_jP = 0$ \qquad $(j\ne 0)$
 \item $\|Pb_0P\| \ge \|b_0\| -\varepsilon $
\end{enumerate} 
\label{lemma:I-free}
\end{lemma} 

\begin{proof} For $x \in X^{\otimes n}$, we define $\length (x) = n$ 
with the convention $\length (a) = 0$ for $a \in A$.  
We write $c$ as a finite sum $c = a + \sum _i S_{x_i}S_{y_i}^*$.
Put $n = 2 \max \{\length (x_i), \length (y_i) ; i\}$.  

\par\noindent
For $j > 0$, each $b_j$ is a finite sum of terms in the form such that 
$$
S_x S_y^* \qquad x \in X^{\otimes (k+j)}, \qquad y \in X^{\otimes k} 
\qquad 0 \le k+j \le n . 
$$
In the case when $j<0$, 
$b_j$ is a finite sum of terms in the form such that 
$$
 S_x S_y^* \qquad x \in X^{\otimes k}, \qquad y \in X^{\otimes (k+|j|)} 
\qquad 0 \le k+|j| \le n . 
$$

We shall identify $b_0$ with an element in 
$L(X^{\otimes n})$.  
Apply Proposition  \ref{prop:free} for 
$T = (b_0)^{1/2}$. Then  
there exists 
$a \in A^+  =C(J)^+ $ with $\|a\|=1$  such that 
\begin{align*}
 \|\phi(a) T\|^2  & \ge \|T\|^2 - \varepsilon,  \\
 a S_xS_y^*a  & = 0 \quad \text{ for any } x \in X^{\otimes i}, 
   \text{ for any }  y \in
 X^{\otimes j}, 0 \le i,j \le N, i \ne j. 
\end{align*}
Define a positive operator $P = a \in A$.  Then 
\[
  \| P b_0 P\| = \| Pb_0^{1/2} \|^2 
                \ge \|b_0^{1/2}\|^2 -\varepsilon 
                = \| b_0 \| -\varepsilon. 
\]
It is evident that $Pb_jP = 0$ for $j \not= 0$.  
\end{proof}

Since we have prepaired technical lemmas addapting to our particular situation, the rest of the proof of our main theorem is a standard one.

\begin{thm}Let $p(z,w)$ be a reduced non-zero polynomial 
in two variables with 
a unique 
factorization into irreducible polynomials:
$
 p(z,w) = g_1(z,w)\cdots g_p(z,w), 
$
where each $g_i(z,w)$ is irreducible and 
 $g_i$ and $g_j \ (i \not= j)$ 
are prime each other.  
We assme that any  $g_i(z,w)$  is not a function only in $z$  or $w$. 
Let $J$ be a $p$-invariant set  with no isolated points.  
Suppose that $p$ is expansive and  free on $J$. 
Then the associated $C^*$-algebra 
${\mathcal O}_p(J)$ 
is simple and purely infinite.
\end{thm}
\begin{proof}
Let $w \in {\mathcal O}_p(J)$ 
be any non-zero positive element.    
We shall show that there exist $z_1$, 
$z_2 \in {\mathcal O}_p(J)$  
such that $z_1^*w z_2 =I$.  
We may assume that $\|w\|=1$.
Let $E : {\mathcal O}_p(J)
 \rightarrow {\mathcal O}_p(J)^{\alpha}$
be the canonical conditional expectation onto the fixed point 
algebra by the gauge action $\alpha$. 
Since $E$ is faithful, $E(w) \not= 0$.  
Choose  $\varepsilon$ such that 
\[
0 < \varepsilon < \frac{\|E(w)\|}{4} \  \text{ and } \ 
\varepsilon (\|E(w)\| -3\varepsilon)^{-1} \leq 1 .  
\]

There exists an element $c \in {\mathcal O}_p(J)^{alg}$  
such that $ \|w - c^* c\| < \varepsilon$ and  
$\|c\| \le 1$.  Let $b = c^*c$. Then $b$ is decomposed  
as a finite sum $b = \sum_j b_j$ with 
$\alpha_t(b_j) =e^{ijt}b_j$.  
Since  $\|b\| \le 1$,  $\|b_0\| = \|E(b)\| \le 1$. 
By Lemma \ref{lemma:I-free}, there exists $P \in A$ with 
$0 \le P \le I$ satisfying $Pb_jP = 0$ \  $(j\ne 0)$
and $\|Pb_0P\| \ge \|b_0\| -\varepsilon $. 
Then we have 

\begin{align*}
\| Pb_0 P \|  & \ge \|b_0 \| -\varepsilon 
              = \|E(b)\| -\varepsilon \\
              &  \ge \|E(w)\| -\|E(w) - E(b)\| -\varepsilon 
              \ge \|E(w)\| -2 \varepsilon . 
\end{align*}
For $T := Pb_0P \in L(X^{\otimes m})$,  
there exists $f \in X^{\otimes m}$ with $\|f\|=1$ such that  
$$
 \|T^{1/2}f \|_2^2 = \|(f|Tf)_A\|  \ge \|T\| -\varepsilon . 
$$
Hence we have 
$\|T^{1/2}f \|_2^2 \ge \|E(w)\| - 3 \varepsilon $.
Define $a = S_f^* T S_f = (f|Tf)_A \in A$.  
Then $\|a\| \ge \|E(w)\| -3 \varepsilon  > \varepsilon$. 
By Lemma \ref{lemma:u-epsilon},  there exists
$n \in \mathbb{N}$ and $u \in X^{\otimes n}$ 
sucn that  
$$
 \| u \| _2 \le (\|a\| - \varepsilon)^{-1/2} \qquad \text{and} 
\quad  S_u^* a S_u = I. 
$$
Then $\|u\| \le (\|E(w)\| -3 \varepsilon)^{-1/2}$. 
We have 
\[
\| S_f^* PwP S_f -a \|  \le \|S_f\|^2 \|P\|^2 \|w -b\| <\varepsilon . 
\]
Therefore  
\[
\|S_u^*S_f^* PwPS_fS_u -I\|  
  < \|u\|^2 \varepsilon 
  \le \varepsilon (\|E(w)\| -3\varepsilon)^{-1} \le 1.
\]
Hence  
$S_u^*S_f^*PwPS_fS_u$ is invertible. Thus there  exists 
$v \in {\mathcal O}_X$ with $S_u^*S_f^*PwPS_fS_u v =I$. 
Put $z_1=S_u^*S_f^*P$ and $z_2=PS_fS_u v$. Then  
$z_1 w z_2 = I$.  

\end{proof}

\noindent
{\bf Remark}.
Schweizer's theorem in \cite{S} also implies that 
${\mathcal O}_p(J)$ is simple. 
Our theorem gives simplicity and pure infiniteness with a 
direct proof. 

$C^*$-algebra 
${\mathcal O}_p(J)$ 
is separable and nuclear, and satisfies the Universal  
Coefficient Theorem. 
Hence the isomorphism class of $C^*$-algebra 
${\mathcal O}_p(J)$ is completely determined by the $K$-group together with the class of the 
unit by the classification theorem by Kirchberg-Phillips \cite{Ki},
\cite{Ph}.

\medskip\par
\noindent
{\bf Example.} Let  $m$ and $n$ be natural numbers. 
Consider $p(z,w)=z^m-w^n$  and $J = {\mathbb T}$.  
If $n$ is not devided by $m$, then 
${\mathcal O}_p(J)$ is simple and purely infinite. 

\medskip\par
\noindent
{\bf Example.} 
Let  $m$ and $n$ be natural numbers and relatively prime. 
Consider  $p(z,w)=(w-z^m)(w-z^n)$ and 
$J = {\mathbb T}$.  Then 
${\mathcal O}_p(J)$ is simple and purely infinite. 

\medskip\par
\noindent
{\bf Example.} Let $R_1(z)=\frac{(z^2 + 1)^2}{4z(z^2 - 1)}$ 
be a rational function given  by Lattes.  Let 
 $R_2(z)=P_2(z)/Q_2(z)$ be any rational function 
with odd degree. 
Consider 
\[
p(z,w) =
 ((4z(z^2 - 1))w-(z^2 + 1)^2)(Q_2(z)w-P_2(z)). 
\]
Let $J=\hat{\mathbb C}$. Then 
${\mathcal O}_p(J)$ is simple and purely infinite. 

\medskip\par
\noindent
{\bf Example.} 
Let $i_1$, $\cdots$, $i_n$, $j_1$, $\cdots$, $j_n$ be natural 
numbers. Assume that $i_k \not= 1$ or $j_k \not= 1$ 
for each $k$.  Suppose that those which are not equal to $1$ 
are relatively prime. Put $J={\mathbb  T}$. Let  
\[
 p(z,w) = (z^{i_1}-w^{j_1})(z^{i_2}-w^{j_2})\cdots (z^{i_n}-w^{j_n}). 
\]
Then 
${\mathcal O}_p(J)$ is simple and purely infinite.

\section{K-groups}

We shall compute K-groups for several examples. 

\medskip\par
\noindent
{\bf Example.} Let $p(z,w) = z^m - w^n $. 
Then $J := {\mathbb T}$ is a $p$-invariant set. 
Consider  the Hilbert bimodule 
$X_p= C(\cC_p)$ over $A = C({\mathbb T})$. 
Then $X_p$ is isomorphic to $A^m$ as a right $A$-module. 
In fact, let $u_i(z,w) = \frac{1}{\sqrt{m}} z^i$ for 
$i = 0,1,\dots, m-1$.  Then $(u_i |u_j)_A = \delta_{i,j}I$ 
and $\{u_0, u_1, \dots, u_{m-1}\}$ is a basis of $X_p$, in the sense that 
$f = \sum_{i=0}^{m-1} u_i(u_i|f)_A$ for any $f \in X_p= C(\cC_p)$. 
Let $a_1(z) = z$ for $z \in {\mathbb T}$. Then 
\[
(\phi (a_1)u_i)(z,w) = \frac{1}{\sqrt{m}} z^{i+1} = u_{i+1}(z,w)
\]
for $i = 0,1,\dots, m-2$. And 
\[
(\phi (a_1)u_{m-1})(z,w) = \frac{1}{\sqrt{m}} z^m  = \frac{1}{\sqrt{m}}w^n 
= (u_0 \centerdot a_1^n)(z,w). 
\]
Therefore, if  we identify $\phi(a) \in L(X_p) = M_n(A)$, then 
$\phi(a_1)_{i,j} = I$ for $i=j+1, j=1,\dots,m-2$, 
$ \phi(a_1)_{0,m-1} = a_1^n$ and  $\phi(a_1)_{i,j} = 0$ for others. 
Let $\phi_1^* : K_1(A) = {\mathbb Z} 
\rightarrow K_1(A) = {\mathbb Z}$.  
Since $[a_1]$ is the generator of  $K_1(A) = {\mathbb Z}$ 
and $\phi_1([a_1])= [a_1^n]$, 
 $\phi_1(k) = nk$ for $k \in {\mathbb Z}$. 

Since  $\phi(I_A) = I_{M_m(A)}$, 
$\phi_0^* : K_0(A) = {\mathbb Z} \rightarrow K_0(A) = {\mathbb Z}$ 
is given by $\phi_1^*(k) = mk$ for $k \in {\mathbb Z}$.    
By a six-term exact sequence due to Pimsner \cite{Pi}, we have 
\[
\begin{CD}
    {\mathbb Z}@>{id - m\cdot}>> {\mathbb Z} 
@>i_*>> K_0({\mathcal O}_p({\mathbb T})) \\
   @A{\delta _1}AA
    @.
     @VV{\delta _0}V \\
   K_1({\mathcal O}_p({\mathbb T})) 
@<<i_*< 0 @<<{id - n\cdot}< {\mathbb Z}
\end{CD}
\]
Therefore 
\par
\noindent
(1)$n = 1$ and $m = 1$: 

$K_0({\mathcal O}_p({\mathbb T})) \cong {\mathbb Z} \oplus {\mathbb Z}$, 
\ \ 
$K_1({\mathcal O}_p({\mathbb T})) \cong {\mathbb Z} \oplus {\mathbb Z}$. 
\par
\noindent
(2)$n = 1$ and $m \not= 1$: 

$ K_0({\mathcal O}_p({\mathbb T})) \cong 
{\mathbb Z} \oplus {\mathbb Z}/(m-1){\mathbb Z}$, \ \ 
$ K_1({\mathcal O}_p({\mathbb T})) \cong {\mathbb Z}$. 
\par
\noindent
(3)$n \not= 1$ and $m = 1$: 

$ K_0({\mathcal O}_p({\mathbb T})) \cong {\mathbb Z}$, \ \ 
$ K_1({\mathcal O}_p({\mathbb T})) \cong 
{\mathbb Z} \oplus {\mathbb Z}/(n-1){\mathbb Z}$. 
\par
\noindent
(4)$n \not= 1$ and $m \not= 1$: 

$ K_0({\mathcal O}_p({\mathbb T})) \cong {\mathbb Z}/(m-1){\mathbb Z}$, \ \ 
$ K_1({\mathcal O}_p({\mathbb T})) \cong {\mathbb Z}/(n-1){\mathbb Z}$. 

\medskip\par
\noindent
{\bf Example.} Let 
$p(z,w) = (w - z^m)(w -  z^n) $ with $(2 \leq m < n)$.  
Then $J := {\mathbb T}$ is a $p$-invariant set. Since 
the set $B(p)$ of branched points in $\cC_p$ is non-empty, 
we need to be 
careful to compute the K-groups $ K_0({\mathcal O}_p({\mathbb T}))$ 
and  $ K_1({\mathcal O}_p({\mathbb T}))$. 

If $z$ is a branched point, then  $w = z^m = z^n$, so that 
$z^{n-m} = 1$. Hence 
$B(p) = \{1,\alpha, \alpha^2, \dots, \alpha^{n-m-1} \}$, 
where $\alpha = e^{2\pi i/(n-m)}$ is a primitive $(n-m)$-th root 
of unity. 
Put 
\[
D(p) = \{(z,w) \in \cC_p \ | \ e(z,w) \geq 2 \} 
     = \{((1,1), (\alpha, \alpha^m), 
\dots, (\alpha^{n-m-1}, \alpha^{m(n-m-1)}) \}
\]  
and any branch index $e(\alpha^k, \alpha^{mk}) = 2$. 
Let $p_1(z,w) = (w - z^m)$ and $p_2(z,w) = (w -  z^n)$. 
Let $X = C(\cC_p)$ and 
\[
Y = C(\cC_{p_1}) \oplus C(\cC_{p_2}) 
  = C(\{(1,z,w) \ | \ p_1(z,w) = 0 \} \cup \{(2,z,w) \ | \ p_2(z,w) = 0 \}).
\]  
We shall embedd $X = C(\cC_p)$ into $Y$ as bimodule over $A = C({\mathbb T})$  
by identifying the points correspoing to the branched points of $p$. 
Let 
\begin{align*}
Z : & = \{f \in Y \ | \ f(1,z,w) = f(2,z,w), | (z,w) \in D(p) \} \\
    & = \{f \in Y \ | \ 
 f(1,\alpha^r, \alpha^mr) = f(2,\alpha^r, \alpha^mr), \  
\ r = 0,1,\dots, n-m-1  \}. 
\end{align*}
Then $Z$ is a closed submodule of $Y$ and we can identify 
 $X$ with $Z$ as bimodule. 

We introduce a basis $\{u_1, \dots, u_{m+n} \}$ of $Y$ as 
follows: 

\par
\noindent
For $i = 1,\dots, m$, 
\[
u_i(1,z,w) = \frac{1}{\sqrt{m}} z^{i-1}, \ \ 
u_i(2,z,w) = 0 ,
\]
and for $i = m+1, \dots, m+n $, 
\[
u_i(1,z,w) = 0, \ \  
u_i(2,z,w) = \frac{1}{\sqrt{n}} z^{i-m-1}. 
\]

Then $(u_i |u_j)_A = \delta_{i,j}I$. Therefore we can 
identify $f \in Y$ with $(f_i)_i \in A^{m+n}$ by 
\[
f_i = (u_i|f)_A, \ \ 
f(k,z,w) = \sum_{i=1}^{m+n} u_i(k,z,w)f_i(w), \ \ k = 1,2.  
\]
We claim that 
$f(1,\alpha^r, \alpha^{mr}) = f(2,\alpha^r, \alpha^{mr})$
if and only if 
\[
\sum_{i=1}^{m} u_i(1,\alpha^r, \alpha^{mr})f_i(\alpha^{mr})
= \sum_{i=m+1}^{m+n} u_i(2,\alpha^r, \alpha^{mr})f_i(\alpha^{mr})
\]
if and only if
\[
\sum_{i=1}^{m} \frac{1}{\sqrt{m}} \alpha^{r(i-1)}f_i(\alpha^{mr})
= \sum_{i=m+1}^{m+n} \frac{1}{\sqrt{n}} \alpha^{r(i-m-1)}f_i(\alpha^{mr}) 
\]
if and only if 
the correspoing vector $(f_1(\alpha^{mr}), \dots, f_{m+n}(\alpha^{mr})) 
\in {\mathbb C}^{m+n}$ is orthgonal to a vector 
\[
n_r :=( \frac{1}{\sqrt{m}}1, \frac{1}{\sqrt{m}}\alpha^r, 
\dots, \frac{1}{\sqrt{m}}\alpha^{r(m-1)},-\frac{1}{\sqrt{n}}1,
-\frac{1}{\sqrt{n}}\alpha^r, \dots, 
-\frac{1}{\sqrt{n}}\alpha^{r(n-1)}) \in {\mathbb C}^{m+n}. 
\]

Let $C := \{\alpha^{mr} \ | \ r = 0,1, \dots, n-m-1 \} 
= \{c_1,c_2, \dots, c_v \}$ and $c_i \not= c_j$, ($i \not= j$). 
For $k = 1,2,\dots, v$, put 
 $C(k) = \{r = 0,1, \dots, n-m-1 \ | \ \alpha^{mr} = c_k \}$. 
If we identify $Y = A^{m+n} = C({\mathbb T})^{m+n}$, then 
\begin{align*}
Z  & =  \{ f = (f_i)_i \in A^{m+n} \ | 
\ (f_i(\alpha^{mr}))_i \text{ is orthogonal to } n_r , 
\text{ for }
\ r = 0,1, \dots, n-m-1   \} \\
&= \cap _{k=1}^v \{f = (f_i)_i \in A^{m+n} \ | 
\ (f_i(c_k))_i \text{ is orthogonal to } n_r 
\text{ in } {\mathbb C}^{m+n} 
\text{ for any }\ r \in B(k) \}.  
\end{align*}
We see that for fixed $k$, the vectors $ n_r \ ( r \in B(k))$ 
are linearly independent. 
Therefore the subspace 
\[
H(k) := \{x= (x_i)_i \in {\mathbb C}^{m+n} \ | \ 
  x \text{ is orthogonal to } n_r, \ r \in B(k) \} 
\]
has dimension 
$m + n - \ ^{\#}C(k) \geq m + n - (n-m) = 2m \geq 2$. 
Let
\[
L(k) := \text{ span } \{ T \in B({\mathbb C}^{m+n}) \ | \ 
T = \theta_{x,y} \text{ for some } x,y \in H(k) \}.  
\]
Therefore  we have an identification 
\[
{\mathcal K}(Z) = \{ f \in C({\mathbb T},
 M_{m+n}({\mathbb C})) \ | \ f(c_k) \in L(k), \ k = 1,\dots, v \}. 
\]
We shall show that the canonical inclusion 
$i : {\mathcal K}(Z) \rightarrow {\mathcal K}(Y) 
\cong M_{m+n}(C({\mathbb T}))$ induces 
the isomorphism 
\[
i_* : K_r({\mathcal K}(Z)) \cong {\mathbb Z} \rightarrow 
        K_r({\mathcal K}(Y)) \cong {\mathbb Z}, \ \ r = 0,1 . 
\]
Let 
\[
J = \{ f \in C({\mathbb T},
 M_{m+n}({\mathbb C})) \ | \ f(c_k)= 0, \ k = 1,\dots, v \}
\] 
and a finite dimensional algebra $Q = \oplus_{k = 1}^v L(k)$. 
Then we have an exact sequence 
\[
0 \rightarrow J \rightarrow {\mathcal K}(Z) 
\overset{\pi}{\rightarrow} Q 
\rightarrow 0.
\]
 Consider the six-term exact sequence
\[
\begin{CD}
  K_0(J)= 0@>>> K_0({\mathcal K}(Z))
@>\pi^*>> K_0(Q)= {\mathbb Z}^v\\
   @A{\delta _1}AA
    @.
     @VV{\delta _0}V \\
   K_1(Q) =0 @<<{\pi}^*< K_1({\mathcal K}(Z)) 
@<<{i^*}< K_1(J) = {\mathbb Z}^{v}
\end{CD}
\] 
For $k = 1,\dots,v$, 
let $q_k \in L(k)$ be a minimal projection and we 
put a projection 
$r_k = (0, \dots,0,q_k,0,\dots,0) \in Q$. Let 
$f_k \in {\mathcal K}(Z)$ a lift of $r_k$ definded 
as a "piecewise linear" map with 
$f_k(c_j) = \delta _{k,j}$. 
Since 
$\delta _0([r_k]) = -[e^{2\pi if_k}]$, 
we obtain that  
\[
\delta _0(n_1,\dots,n_v) 
= (n_v-n_1, n_1-n_2,n_2 -n_3,\dots, n_{v-1}-n_v). 
\]
Since ${\rm Im} \pi^* =  \Ker \delta_0 
= \{(n,\dots,n) \in {\mathbb Z}^r \ | \ n \in {\mathbb Z} \}
\cong {\mathbb Z}$ and $\pi^*$ is one to one, we see 
$\pi^* : K_0({\mathcal K}(Z)) \cong {\mathbb Z} 
\rightarrow {\mathbb Z}^r$ is given by 
$\pi^*(n) = (n,\dots,n)$. Since ${\rm Im} \  \delta_0 
= \Ker  \ i^*$ 
and $i^*$ is onto,  
$i^* : K_1(J) \cong {\mathbb Z}^r \rightarrow 
K_1({\mathcal K}(Z)) \cong {\mathbb Z}$ is given by 
$i^*(n_1,\dots,n_v) = n_1 + \dots + n_v$.  
Let $p \in  C({\mathbb T}, M_{m+n}({\mathbb C}))$ be 
a projection such that $p(t)$ is a rank one projection for any 
$t \in {\mathbb T}$ and $p(c_k) \in L(k)$ for $k = 1,\dots, v$. 
Then $[p]$ is a generator of $K_0({\mathcal K}(Z)) 
\cong {\mathbb Z}$ and 
also a generator of $K_0({\mathcal K}(Y)) \cong {\mathbb Z}$. Let 
$c_k = e^{2\pi i \theta_k}$ with 
$0 \leq \theta_1 \leq \dots \leq  \theta_v$. 
Let $u \in C({\mathbb T}, M_{m+n}({\mathbb C}))$ be  
a unitary such that $u(e^{2\pi i t}) = e^{2\pi i t/{\theta_1}}$ 
for $0 \leq t \leq {\theta_1}$ and $u(e^{2\pi i t}) = 1$ for 
others. Then $[u]$ is a 
generator of $K_1({\mathcal K}(Z)) \cong {\mathbb Z}$ and 
also a generator of $K_1({\mathcal K}(Y))$.
Therefore we conclude  that 
$i_* : K_r({\mathcal K}(Z)) \cong {\mathbb Z} \rightarrow 
        K_r({\mathcal K}(Y)) \cong {\mathbb Z}, \ \ r = 0,1 
$ is an isomorphism.  

Since $I_X = \{ f \in C({\mathbb T}) \ | \ f(\alpha^k) = 0 
\text{ for } k = 0,1, m-n-1 \}$, we have $K_0(I_X) = 0$ and 
$K_1(I_X) = {\mathbb Z}^{m-n}$. 

Therefore we can identify the left action 
$\phi: I_X \rightarrow K(X) = K(Z)$ with 
$\phi_1 \oplus \phi_2: I_X \rightarrow 
K(Y) = K(C(\cC_{p_1}) \oplus C(\cC_{p_2}))$ 
on the level of K-groups.
Hence $\phi^*: K_1(I_X) = {\mathbb Z}^{m-n}  
\rightarrow K_1(A)= {\mathbb Z}$ is given by 
$\phi^*(x_1,\dots,x_{m-n}) = \sum_{i = 1}^{m-n} 2x_i$.  

By a six-term exact sequence, we have 
\[
\begin{CD}
   K_0(I_X)= 0@>{j^* - \phi^*}>> K_0(A) = {\mathbb Z} 
@>i_*>> K_0({\mathcal O}_p({\mathbb T})) \\
   @A{\delta _1}AA
    @.
     @VV{\delta _0}V \\
   K_1({\mathcal O}_p({\mathbb T})) 
@<<i_*< K_1(A) = {\mathbb Z} @<<{j^* - \phi^*}< 
K_1(I_X) = {\mathbb Z}^{m-n}
\end{CD}
\]
The canonical inclusion map $j: I_X 
\rightarrow A = C({\mathbb T})$ induces 
$j^*: K_1(I_X) = {\mathbb Z}^{m-n} \rightarrow K_1(A) = {\mathbb Z}$ 
with $j^*(x_1,\dots,x_{m-n}) = \sum_{i = 1}^{m-n} x_i$.  
Therefore we have 
\[
 K_0({\mathcal O}_p({\mathbb T})) = {\mathbb Z}^{m-n}, \ \ 
\text{ and } \ \    K_1({\mathcal O}_p({\mathbb T})) = 0.
\]

\medskip\par
\noindent
{\bf Example.} Let 
$p(z,w) = (w - z^{m_1})(w -  z^{m_2})\dots (w - z^{m_r}) $ 
and $m_1, \dots, m_r$ are all different, where $r$ is the 
number of irreducible components. 
Then $J := {\mathbb T}$ is a $p$-invariant set. 
Let $b = \ ^{\#}B(p)$ be the number of the branched poins. 
By a similar 
calculation, we have 
\[
 K_0({\mathcal O}_p({\mathbb T})) = {\mathbb Z}^{b}, \ \ 
\text{ and } \ \   
 K_1({\mathcal O}_p({\mathbb T})) = {\mathbb Z}/(r-1){\mathbb Z}.
\]

\end{document}